\documentclass[twocolumn]{autart} 
\usepackage{graphicx}

\usepackage{epstopdf}
\usepackage{amsmath,amssymb,amsfonts,bm,amscd}
\usepackage{epsfig,subfigure}
\usepackage{verbatim}
\usepackage{color}

\newtheorem{lemma}{Lemma}

\newtheorem{theorem}{Theorem}

\newtheorem{example}{Example}

\begin{document}

\begin{frontmatter}

\title{Infinite Chains of Kinematic Points\thanksref{footnoteinfo}}

\thanks[footnoteinfo]{This paper was not presented at any IFAC 
meeting. The corresponding  author is Bruce Francis.}

\author[Abie]{Avraham Feintuch}\ead{abie@math.bgu.ac.il},    
\author[Bruce]{Bruce Francis}\ead{bruce.francis@utoronto.ca}       

\address[Abie]{Department of Mathematics, Ben-Gurion University of the Negev, Israel}  
\address[Bruce]{Department of Electrical
and Computer Engineering, University of Toronto, Canada}             

\maketitle


\begin{abstract}
In formulating the stability problem for  an infinite chain of cars, state space is traditionally taken to be  the Hilbert space $\ell^2$, wherein the displacements of  cars from their equilibria, or the velocities from their equilibria,  are taken to be square summable. But this obliges the displacements or velocity perturbations of cars that are far down the chain to be vanishingly small and leads to anomalous behaviour. In this paper an alternative formulation is proposed wherein state space is the Banach space $\ell^\infty$, allowing the displacements or velocity perturbations of  cars from their equilibria to be merely bounded.
\end{abstract}
  
  \end{frontmatter}

\section{Introduction}\label{section_intro}

In studying the formation of a very large number of vehicles, one approach is instead to model an infinite number of vehicles \cite{BamPagDah02}, \cite{DanDul03},  \cite{MelKuo71},
 \cite{MotJad08}. The question then arises as to what mathematical framework to take so that the latter model correctly describes the behaviour of the  former. The purpose of this paper is to suggest that the Hilbert space framework usually adopted is not always appropriate and to suggest an alternative.

\begin{figure}[htbp] 
   \begin{center}
   \includegraphics[width=3.2in]{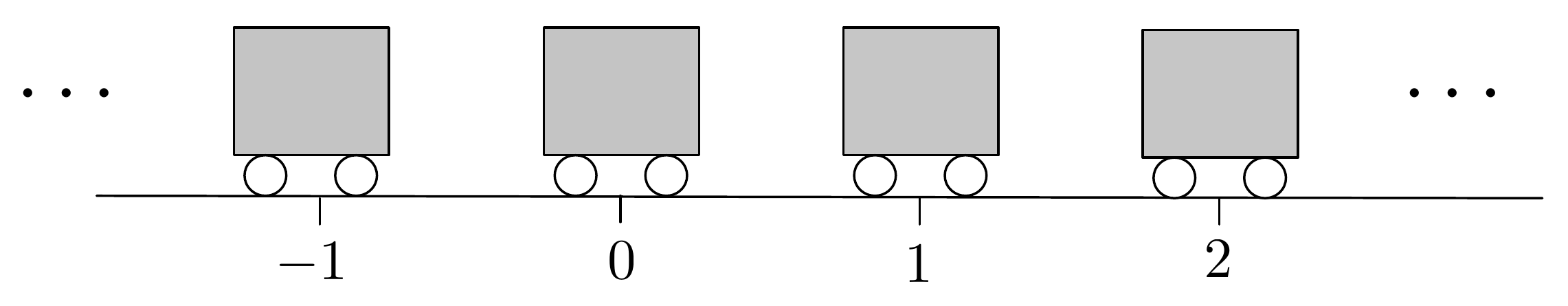} 
     \caption{Infinite chain of cars}
   \label{fig_cars}
     \end{center}
\end{figure}

Consider the infinite chain of cars in Figure \ref{fig_cars}.
The cars are modelled as points on the real line $\mathbb{R}$ and are numbered by the integers.  The position of car $n$ is denoted by $q_n \in \mathbb{R}$. 
We take the simplest model of a car, a kinematic point:
\begin{equation*}
\dot{q}_n=u_n, \ \ \ n \in \mathbb{Z}.
\end{equation*}
With
nearest-neighbour interaction, the control velocity  would be of the form
\begin{equation*}
u_n=f(q_{n+1}-q_n,q_{n-1}-q_{n}), \ \ \ n \in \mathbb{Z},
\end{equation*}
where typically the function $f$ is linear and the same for all $n$. 
Thus
\[
\dot{q}_n=f(q_{n+1}-q_n,q_{n-1}-q_{n}).
\]
The cars are nominally spaced a unit distance  apart. It is assumed that  $q_n=n$ is an equilibrium of the system, that is,
$f(1,-1)=0.$
 Let $p_n$ denote the displacement of car $n$ away from its equilibrium position: $p_n=q_n-n$. Thus the nominal displacements are $p_n=0$. With $f$ linear it follows that $p_n$ satisfies the same equation as $q_n$:
\begin{equation}
\dot{p}_n=f(p_{n+1}-p_n,p_{n-1}-p_{n}), \ \ \ n \in \mathbb{Z},
\label{eq_mass_spring}
\end{equation}
The difference between the two models is that it is natural to take bounded initial conditions in the $p$-model. 
Thus the model is an infinite number of coupled differential equations. 

Let $p(t)$ denote the infinite vector of displacements at time $t$, that is, the components of $p(t)$ are $p_n(t)$, $n \in \mathbb{Z}$. In this paper we are  interested in the question of stability---does $p(t)$ converge as $t \rightarrow \infty$ and if so in what sense?
In the existing literature, e.g., \cite{BamPagDah02}, \cite{DanDul03}, the state space for $p(t)$ is the Hilbert space $\ell^2$ of square-summable sequences, the advantage of this setting being that Fourier transforms can be exploited. But this assumption requires that $p_n(t) \rightarrow 0$ as $n$ goes to $\pm \infty$, for every $t$. This seems to be an unjustified assumption to make at the start of a stability theory, before anything has been proved: If we want to know about the behaviour of $p(t)$ as $t\rightarrow \infty$ there is no justification in limiting $p(0)$ to satisfy  $p_n(0) \rightarrow 0$ as $n \rightarrow \pm \infty$. Therefore we take the state space to be  the Banach space $\ell^\infty$ of bounded sequences. Then $p_n(0)$ can all be of roughly  equal magnitude, or they can be randomly distributed in an interval, etc. The only requirement is that $p_n(0)$ lie in some interval independent of $n$. The goal of this paper is to develop a stability theory in this context,  an $\ell^\infty$ theory, and to show that it is different from the $\ell^2$ theory. We illustrate with two examples.

\begin{example}
{\rm Suppose each car heads toward the sum of the relative displacements to its two neighbours: 
\begin{align*}
\dot{q}_n  &=  (q_{n+1}-q_n) +(q_{n-1}-q_{n})\\
&=  q_{n+1}+q_{n-1}-2q_n.
\end{align*}
It follows that $p_n$ satisfies the same equation:
\[
\dot{p}_n  =  p_{n+1}+p_{n-1}-2p_n.
\]
 Let $p$ denote the infinite vector of displacements.
Thus $p=0$ is an equilibrium. If $p(0) \in \ell^2$, it turns out (proved in the paper) that $p(t)$  converges to zero, that is, the cars return to their original positions! 
  But why should the infinite chain  behave in this way? After all, the cars are not fitted with global sensors to know where the origin is. This anomaly is caused entirely by taking $p(0)$ in $\ell^2$.}
  \end{example}
  
\begin{example}
{\rm  Consider again an infinite chain of cars, but where the control objective is to maintain a constant distance between
the cars and a constant velocity for each car. Let the desired distance  between two consecutive cars be $d$ and
 the desired velocity of each be $v_d$.  Suppose that for $t<0$  the cars are  spaced exactly distance $d$ apart and are all moving in the same direction at speed $v_d$. That is an equilibrium situation. Let the velocities of the cars be denoted $v_n(t)$, where $n$ ranges over all integers. Now suppose that at time $t=0$ every car suddenly speeds up by $1\%$. This is a perturbation away from the equilibrium. At $t=0$ the velocity of every car is $v_d+0.01v_d$. The perturbations $\hat{v}_n(0)=0.01 v_d$ are not  square-summable, that is,
\[
\sum_{n=-\infty}^\infty \hat{v}_n(0) ^2 = \infty.
\]
Therefore the vector of perturbations $\hat{v}(0)$ is not in $\ell^2$ but rather is in $\ell^\infty$. That is the situation we are discussing in the paper and that is not handled by the $\ell^2$ formulation.}
\end{example}

The literature review 
is postponed until the end of the paper, where we have suitable notation.

\section{Mathematical Preliminaries}

The signals that we deal with are denoted, for example, by $x(t)$, where $t$ denotes time and $x$ is a vector with an infinite number of components, $x_n$, $n \in \mathbb{Z}$. The meaning is that $x_n(t)$ is the state vector of car $n$. For simplicity, the dimension of $x_n(t)$ is just 1. Thus for each $t$, $x(t)$ is the state vector of the entire chain. 

\subsection{Spaces, Operators, and Spectra}

In this subsection, $x$, $y$, etc.\  will denote generic sequences of real numbers with components $x_n$, $y_n$, etc.
Let $\mathbb{R}^\infty$ denote the space of such sequences. The index $n$ runs over the set of all integers. The space $\mathbb{R}^\infty$ is an infinite-dimensional vector space. It does not have a norm, though there is a natural topology arising from componentwise convergence. The Hilbert space $\ell^2$ of square-summable sequences and the Banach space $\ell^\infty$ of bounded sequences are both subspaces of $\mathbb{R}^\infty$. The space $\ell^2$ is based on the inner product
\[
\langle x,y \rangle =\sum_{-\infty}^{\infty} x_ny_n.
\]
The induced norm is
\[
\| x \|_2 = \left( \sum_n x_n^2 \right)^{1/2}.
\]
And $\ell^\infty$ is based on the norm
\[
\| x \|_\infty = \sup_n |x_n|.
\]
Of course every square-summable sequence is bounded and therefore $\ell^2$ is a subset of $\ell^\infty$.

Let $\mathcal{X}$ be a Banach space. For us it will be either $\ell^2$ or $\ell^\infty$. The space of bounded linear operators on $\mathcal{X}$ is denoted $\mathcal{B}(\mathcal{X})$, or just $\mathcal{B}$.  Let $A \in \mathcal{B}$. A complex number $\lambda$ is a {\bf regular point} of $A$ if  $(\lambda I-A)^{-1} \in \mathcal{B}$. The set of regular points is the {\bf resolvent set}, and its complement, $\sigma(A)$, the {\bf spectrum} of $A$. An eigenvalue is, as for matrices, a complex number $\lambda$ for which there exists a nonzero $x$ in $\mathcal{X}$ such that $Ax=\lambda x$. Eigenvalues, if they exist, certainly belong to the spectrum, but  $\lambda I-A$ may fail to have a bounded inverse for other complex numbers $\lambda$ than eigenvalues. The spectrum is always nonempty, closed, and  contained in the disk  $| z | \leq r_A$, where the {\bf spectral radius} $r_A$ is given by
\[
r_A={\limsup}_{n \rightarrow \infty} \| A^n \|^{1/n}.
\]

\subsection{The Bilateral Right Shift}

The {\bf bilateral right shift} $U$ is the linear transformation on  $\mathbb{R}^\infty$ defined by 
\[
y=Ux, \ \ y_n=x_{n-1}.
\]
Its inverse is the left shift:
\[
y=U^{-1}x, \ \ y_n=x_{n+1}.
\]
When restricted either to $\ell^2$ or to $\ell^\infty$, $U$ is a bounded operator, that is, $U$ belongs to both $\mathcal{B}(\ell^2)$ and $\mathcal{B}(\ell^\infty)$. Its properties in these two spaces are somewhat different.

\begin{lemma}\label{lemmaU}
(Properties of $U$)
\begin{enumerate}
\item
Consider $U$ as an operator in $\mathcal{B}(\ell^2)$. Its spectrum equals the unit circle, but $U$ has no eigenvalues. The kernel (nullspace) of $U-I$ equals the zero subspace.
\item
Consider $U$ as an operator in $\mathcal{B}(\ell^\infty)$. Its spectrum equals the unit circle and every point in the spectrum is an eigenvalue.  The kernel of $U-I$ is the 1-dimensional subspace spanned by the vector $\bf{1}$, whose components are all 1.
\end{enumerate}
\end{lemma}

\noindent
{\bf Proof} \ (1)
This result is standard, e.g., \cite{Hal82}. As an operator on $\ell^2$, $U$ has no eigenvalues. To see this, suppose $Ux=\lambda x$, $x \in \ell^2$, $x \neq 0$. Then
\[
x_{n-1}=\lambda x_n, \ \ \ n \in \mathbb{Z}.
\]
Without loss of generality starting with $x_0=1$, we have by iterating backward in the index that $x_{-m}=\lambda^m$, $m>1$. Since $x \in \ell^2$, so $\lim_{m \rightarrow \infty} \lambda^m =0$ and so $| \lambda |<1$. But by iterating forward in the index from $x_0=1$ we conclude that $|\lambda|>1$. This inconsistency shows there is no $\lambda$. A nonzero vector in the kernel of $U-I$ would be an eigenvector and 1 would be an eigenvalue. But there are no eigenvalues.

(2) Since $\| U^n \|=1$ for all $n$, the spectral radius is
\[
r_U=\lim_n \| U^n \|^{1/n} = 1.
\]
Thus $\sigma(U)$ is contained in the closed unit disk. In fact, 1 is an eigenvalue, with eigenvector  all 1's.
Similarly, for every real $\theta$, $\mathrm{e}^{j\theta}$ is an eigenvalue, with eigenvector $x=(x_n)$ defined by
\[
x_0=1, \ \ \ x_{n-1}=\mathrm{e}^{j\theta}x_n.
\]
So the unit circle is contained in $\sigma(U)$.  The norm of $U^{-n}$ also equals 1 for every $n$, and therefore the spectral radius of $U^{-1}$ equals 1 too. The spectra of $U$ and $U^{-1}$ are reciprocals. This fact and the equalities $r_U=r_{U^{-1}}=1$ show that $\sigma(U)$ has no points in $|z|<1$. Thus $\sigma(U)$ and $\sigma(U^{-1})$ both equal the unit circle. Finally, the kernel of $U-I$ is the eigenspace for the eigenvalue 1. It is easy to check that all eigenvectors are constant. 
This concludes the proof.

\bigskip
Let $A \in \mathcal{B}(\ell^2) \cap \mathcal{B}(\ell^\infty)$. The spectrum of $A \in \mathcal{B}(\ell^2)$ and 
 the spectrum of $A \in \mathcal{B}(\ell^\infty)$ are not always equal (see Example~\ref{example_non_inv}) but,
as we saw, they {\bf are} for $A=U$. And this extends to the case where $A$ is a polynomial in $U,U^{-1}$, for example
\[
A=a_2U^2+a_1U+a_0I+a_{-1}U^{-1}+a_{-2}U^{-2}.
\]
The spectrum of this $A$ equals 
\[
\{ a_2z^2+a_1z+a_0+a_{-1}z^{-1}+a_{-2}z^{-2} : |z|=1\}.
\]

\subsection{Differential Equations in Banach Space}

Here we summarize the theory of Dalecki\u{\i} and Kre\u{\i}n \cite{DalKre70}. The results are for a general Banach space $\mathcal{X}$. 
We shall need the results for $\mathcal{X}=\ell^\infty$ and $\ell^2$ (a Hilbert space is also a Banach space).

The exponential $\mathrm{e}^{A}$ can  be defined by the series
\[
\mathrm{e}^{A}=I+A+\frac{1}{2!}A^2+\cdots.
\]
The operator $\mathrm{e}^{A}$ belongs to $\mathcal{B}$ whenever $A$ does.
It follows that the function $t \mapsto \mathrm{e}^{At}$ is differentiable and satisfies
\[
\frac{d}{dt} \mathrm{e}^{At}=A\mathrm{e}^{At}=\mathrm{e}^{At}A.
\]
It also satisfies $\left. \mathrm{e}^{At}\right|_{t=0}=I$. Consequently for any $x_0 \in \mathcal{X}$, $x(t)=\mathrm{e}^{At}x_0$ satisfies the equation
\[
\dot{x}=Ax, \ \ \ x(0)=x_0.
\]
In fact, it is the unique solution among differentiable functions.

 The spectral mapping theorem holds in this general context. Let $A \in \mathcal{B}$. Let $K_A$  denote the class of functions $\phi(z)$ that are piecewise analytic on $\sigma(A)$. This means that 1) the domain of definition of $\phi$ consists of a finite number of open connected components whose union contains $\sigma(A)$, each component containing at least one point of $\sigma(A)$; and 2) The function $\phi$ is analytic in each component of its domain of definition.
Then
\[
\sigma(\phi(A))=\phi(\sigma(A)),
\]
which says that the spectrum of the operator $\phi(A)$ equals the set of points  $\phi(z)$ as $z$ ranges over the spectrum of $A$.
In particular, the spectrum of $\mathrm{e}^A$ equals the set of points  $\mathrm{e}^z$ as $z$ ranges over the spectrum of $A$.

Theorem 4.1, page 26, of \cite{DalKre70} provides the following key fact for studying stability.

\begin{theorem}
Let $A \in \mathcal{B}$. If $\sigma(A)$ lies in the open half-plane $\mathrm{Re} \ \lambda < a$, then there exists a constant $b$ such that 
\[
\| \mathrm{e}^{At} \| \leq b \mathrm{e}^{a t}, \ \ \ t \geq 0.
\]
Conversely, if such $b$ exists, then $\sigma(A)$ lies in the closed half-plane $\mathrm{Re} \ \lambda \leq a$.
\end{theorem}

Thus, just as for matrices, if $\sigma(A)$ lies in the open half-plane $\mathrm{Re} \ \lambda < a$ and $a$ is negative, then $\mathrm{e}^{At}$ converges to 0 as $t \rightarrow \infty$,

\subsection{Spatial Invariance and the $\ell^2$-Induced Norm}

Recall that the Fourier transform of $x \in \ell^2$ is the function
\[
X(\mathrm{e}^{j \omega})=\sum_n x_n \mathrm{e}^{-j \omega n}.
\]
This function belongs to the Hilbert space, denoted by $\mathcal{L}^2(S^1)$, of square-integrable functions on the unit circle, $S^1$ (the notation suggests the 1-dimensional unit sphere). The mapping 
\[
F: \ell^2 \longrightarrow \mathcal{L}^2(S^1), \ \ \ F: x \mapsto X
\]
is the {\bf Fourier operator}. It is an isomorphism of Hilbert spaces.

Consider an infinite chain modelled by the first-order equation 
\[
\dot{x}(t)=Ax(t).
\]
The derivative is with respect to time $t$ and $A$ is an operator in $\mathcal{B}(\ell^2)$. The infinite chain is said to be {\bf spatially invariant} if $A$ commutes with $U$, i.e., $A$ is a Toeplitz operator. Then $FAF^{-1}$ is the operator on $\mathcal{L}^2(S^1)$ of multiplication by a function $G(\mathrm{e}^{j\omega})$, the spatial transfer function. Likewise, $F\mathrm{e}^{At}F^{-1}$ is the operator on $\mathcal{L}^2(S^1)$ of multiplication by  $\mathrm{e}^{G(\mathrm{e}^{j \omega})t}$.
Then for each fixed $t$, the $\ell^2$-induced norm $\| \mathrm{e}^{At} \|$ is given by
\begin{align}
\| \mathrm{e}^{At} \| & =  \max_{\omega} \left| \mathrm{e}^{G(\mathrm{e}^{j \omega})t} \right|  \nonumber \\ 
& =  \max_{\omega}  \mathrm{e}^{\mathrm{Re} \ G( \mathrm{e}^{j \omega})t}  \nonumber \\  
& =    \mathrm{e}^{\max_{\omega}  \mathrm{Re} \ G( \mathrm{e}^{j \omega})t} .  \label{eq_norm} \\ \nonumber
\end{align}

\section{Serial Pursuit and Rendezvous}

In this section we take the non-symmetric coupling where car $n$ pursues car $n-1$, for every $n$, according to
\begin{equation}
\dot{q}_n=q_{n-1}-q_n.
\label{eq_serial}
\end{equation}
This setup is not quite like the one in the introduction with symmetry, because  $q_n(0)=n$ is not an equilibrium in (\ref{eq_serial}). So we shall not assume the cars are initially spread out to infinity but rather are all within a bounded interval: $| q_n(0)| \leq B$ for some $B$ and all $n$. We are interested in whether the cars {\bf rendezvous}, that is, all converge to the same location. 

The vector form of (\ref{eq_serial}) is
\begin{equation}
\dot{q}=(U-I)q, \label{pursuit}
\end{equation}
whose solution is $q(t)=\mathrm{e}^{(U-I)t} q(0)$. This operator $U-I$ is represented by the infinite Toeplitz matrix
\[
\left[ \begin{array}{rrr|rrrr}
 & \vdots & \vdots & \vdots & \vdots & \vdots &  \\ 
\dots & -1 & 0 & 0 & 0 & 0 & \dots \\
\dots & 1 & -1 & 0 & 0 & 0 & \dots \\ \hline
\dots & 0 & 1 & -1 & 0 & 0 & \dots \\
\dots & 0 & 0 & 1 & -1 & 0 & \dots \\
\dots & 0 & 0 & 0 & 1 & -1 & \dots \\
 & \vdots & \vdots & \vdots & \vdots & \vdots & 
\end{array} \right] .
\]
The vertical and horizontal lines in the matrix separate the index range $n<0$ from the range $n \geq 0$.
The spectrum of this operator is the circle of radius 1, centre $-1$. Thus 0 is in the spectrum. Perhaps contrary to one's intuition, $-1$ is not in the spectrum.

 It's enlightening to compare the infinite chain with a finite one, as we do  in three examples.
 
\begin{example} \label{example_20}
{\rm Suppose  there are only finitely many cars, in fact, only three cars: $n = 0,1,2$. With only three cars we need a boundary condition for $n=0$ because there is no $n=-1$. One possibility is that car 0 is tethered and hence stationary:
\[
\dot{q}=Aq, \ \ \  A=\left[ \begin{array}{rrr}
0 & 0 & 0 \\
1 & -1 & 0 \\
0 & 1 & -1
\end{array} \right] .
\]
Cars 1 and 2 converge to the stationary car 0.}
\end{example}

\begin{example} \label{example_21}
{\rm Continuing with the same setup, take the boundary condition to be that car 0 can see the global origin and heads for it: 
\[
\dot{q}=Aq, \ \ \  A=\left[ \begin{array}{rrr}
-1 & 0 & 0 \\
1 & -1 & 0 \\
0 & 1 & -1
\end{array} \right] .
\]
All cars converge to the origin.}
\end{example}

\begin{example} \label{example_22}
{\rm Continuing still with the same setup, take the boundary condition to be that car 0 can see car 2 and heads for it: 
\[
\dot{q}=Aq, \ \ \  A=\left[ \begin{array}{rrr}
-1 & 0 & 1 \\
1 & -1 & 0 \\
0 & 1 & -1
\end{array} \right] .
\]
All cars converge to the average of their starting points.}
\end{example}

We return now to the subject of the paper---infinitely many cars with no boundary condition. 
We shall see that if $q(0)$ belongs to $\ell^2$, then $q(t)$ converges to 0, just as in Example~\ref{example_21}. 
This behaviour, where cars without global sensing capability rendezvous at the origin of the global coordinate system, is an anomaly caused by the Hilbert space hypothesis.

\begin{theorem}\label{theorem00}
With reference to (\ref{pursuit}),
the $\ell^2$-induced norm of  $\mathrm{e}^{(U-I)t}$ satisfies $\| \mathrm{e}^{(U-I)t}\| = 1$ for all $t \geq 0$.
For every $q(0) \in \ell^2$, the $\ell^\infty$-norm of $q(t)$ converges to 0 as $t$ tends to $\infty$; 
in addition, $q(t)$ converges  to zero weakly, that is, the $\ell^2$ inner product $\langle q(t), y \rangle$ 
converges to zero as $t \rightarrow \infty$ for every $y$ in $ \ell^2$.
\end{theorem}

The proof is a minor modification of the proof of Theorem~\ref{theorem1} to follow, and hence is omitted.

\bigskip
The case where  $q(0)$ instead belongs to $\ell^\infty$ is significantly more interesting. In this case  all the points start merely
in some ball centred 
at the origin, i.e., $|q_n(0)| \leq \| q(0) \|_\infty$.  We don't have a complete theory on this $\ell^\infty$ problem; what we do have are six results presented in the subsections to follow.

\subsection{Rendezvous with weakened initial conditions}

Our first result relaxes the assumption $q(0) \in \ell^2$ to merely $\lim_{n \rightarrow \pm \infty}q_n(0)=0$. We will see that the cars again rendezvous at the origin. By linearity, this is the equivalent to saying that if $\lim_{n \rightarrow \pm \infty}q_n(0)=c$, then the cars rendezvous at the location $c$.
 

\begin{lemma}\label{theorem12}
Assume $q(0) \in \ell^{\infty}$.
If $q_n(0)$ tends to 0 as $n$ tends to $\pm \infty$, then
$q(t)$ converges in $\ell^\infty$ to $0$ as $t \rightarrow \infty$.
\end{lemma}

\noindent
{\bf Proof} \
Following signal processing notation, let $\delta$ denote the unit impulse in $\ell^\infty$; that is, $\delta_n=0$ for all $n$ except that $\delta_0=1$. Consider the case where $q(0)=\delta$. We have
 \begin{align*}
 \| q(t) \|_\infty &= \| \mathrm{e}^{(U-I)t}\delta\|_\infty \\
 &= \mathrm{e}^{-t}\| \mathrm{e}^{Ut} \delta \|_\infty\\
 &= \mathrm{e}^{-t}\| \delta +tU \delta +(t^2/2!) U^2 \delta + \cdots \|_\infty \\ 
 &= \mathrm{e}^{-t}\sup_{k \geq 0}\left|  \frac{t^{k}}{k!} \right| \\
 &  \rightarrow 0 \mbox{ as } t\rightarrow \infty.
 \end{align*}
 
Next,  if $q(0)$ is a finite linear combination of coordinate vectors, that is, of $\{ U^k \delta \}$, it follows by linearity and the triangle
inequality that 
$\| q(t) \|_\infty \rightarrow 0$ as $t\rightarrow \infty$.
The closed linear span of finite linear combinations of the coordinate vectors in $\ell^{\infty}$ is the subspace $c_{0}$
of vectors $f$ such that $f_n\rightarrow 0$ as $ n\rightarrow \pm \infty$. 
Note that the semigroup $\{\mathrm{e}^{Ut}\}_{t \geq 0}$ satisfies, for $f\in \ell^\infty$,
\begin{align*}
\| \mathrm{e}^{Ut}f \|_\infty &= \| f +tUf +(t^2/2!)U^2f + \cdots \|_\infty \\
& \leq \| f\|_\infty + t \| Uf\|_\infty +(t^2/2!)\| U^2f\|_\infty + \cdots  \\
& = \mathrm{e}^{t} \| f \|_\infty .
\end{align*}
And so $\| \mathrm{e}^{Ut} \| \leq  \mathrm{e}^{t}$.
Thus $\| \mathrm{e}^{(U-I)t}\|  = \mathrm{e}^{-t}\| \mathrm{e}^{Ut}\| \leq 1$. 
Since $c_0$ is invariant under $U-I$, and therefore under $\mathrm{e}^{(U-I)t}$,
we can  apply the uniform boundedness principle to obtain that for any $q(0)\in c_0$, $q(t) \rightarrow 0$ 
as $t\rightarrow \infty$. This completes the proof.

\subsection{Convergence of car $n$ implies that of car $n+1$}

We will see later that a car's position doesn't necessarily converge, as $t \rightarrow \infty$, under just the assumption that $q(0) \in \ell^\infty$.
Here we show that if car $n$ converges to some location, then cars $n,n+1, n+2,\dots$ rendezvous.  This is almost obvious, because car $n_0+1$ is pursuing car $n_0$, car $n_0+2$ is pursuing car $n_0+1$, and so on. 

\begin{lemma}
If the limit $\lim_{t \rightarrow \infty} q_n(t)$ exists for $n=n_0$, then it exists for every $n>n_0$ and all the limits are equal.
\end{lemma}

\noindent
{\bf Proof} \  A direct computation shows that the
$n^{\rm th}$ coordinate of $q(t)$ is given by 
\begin{equation}
q_n(t)=\mathrm{e}^{-t} \sum_{k=0}^{\infty}q_{n-k}(0)\frac{t^{k}}{k!}.
\label{eq_qn}
\end{equation}
Alternatively
\[
q_n(t)=\frac{h_n(t)}{\mathrm{e}^{t}}, \ \ \ h_n(t)= \sum_{k=0}^{\infty}q_{n-k}(0)\frac{t^{k}}{k!}.
\]
Notice that $\dot{h}_{n+1}=h_n$. Therefore if $\lim_{t \rightarrow \infty} q_n(t)=c$, then\footnote{
To apply l'H\^{o}pital's rule in the equations to follow, we need $h_{n+1}(t) \rightarrow \infty$ as $t \rightarrow \infty$. If this isn't the case, simply perturb $q_n(0)$ to $q_n(0)+\varepsilon$ for  a small positive $\varepsilon$, that is, translate all the points. Then $h_{n+1}(t)$ is perturbed to
$h_{n+1}(t)+\varepsilon \mathrm{e}^t$.}
\begin{align*}
\lim_{t \rightarrow \infty} q_{n+1}(t) &= \lim_{t \rightarrow \infty}\frac{h_{n+1}(t)}{\mathrm{e}^{t}} \\
&=\lim_{t \rightarrow \infty} \frac{\dot{h}_{n+1}(t)}{\mathrm{e}^{t}} \ \ \ \mbox{by l'H\^{o}pital's rule}\\
&=\lim_{t \rightarrow \infty} \frac{{h}_{n}(t)}{\mathrm{e}^{t}}\\
&=\lim_{t \rightarrow \infty} q_{n}(t) \\
& = c.
\end{align*}
This concludes the proof.

We don't know if the result extends to $n<n_0$.

\subsection{A sufficient condition for rendezvous}

For every $t>0$, the series (\ref{eq_qn}) converges, that is, $q_n(t)$ is well-defined. However the limit $\lim_{t \rightarrow \infty} q_n(t)$ may or may not exist, depending on $q(0)$. 
Here we give  an example where rendezvous occurs without the strong assumption on $q(0)$ that is in  Lemma~\ref{theorem12}.

\begin{example}
{\rm Take
\begin{align*}
q(0)&=(\dots ,q_{-2}(0),q_{-1}(0) | q_0(0),q_1(0),\dots)\\
&= (\dots , 0,0,0 | 1,0,-1,0,1,0,-1,\dots),
\end{align*}
for which
\[
q_n(t)= \mathrm{Re} \ \mathrm{e}^{(j-1)t}.
\]
Thus $q_n(t)$ converges to 0 as $t \rightarrow \infty$, i.e., the cars rendezvous at the origin. Along a similar line, let $a$ be a real number that is not a rational multiple of $\pi$. Kronecker's density theorem says that the sequence $\{ \mathrm{e}^{jna} \}_{n \in \mathbb{Z}}$ is dense on the unit circle, that is, every point on the circle is an accumulation point for the sequence. For $q_n(0)=\mathrm{e}^{jna}$ we have
\begin{align*}
\lim_{t \rightarrow \infty} q_n(t) &= \lim_{t \rightarrow \infty}\mathrm{e}^{-t} \sum_{k=0}^{\infty} \mathrm{e}^{j(n-k)a} \frac{t^k}{k!} \\
&= \lim_{t \rightarrow \infty} \mathrm{e}^{-t}  \mathrm{e}^{jna}\mathrm{e}^{\left(\mathrm{e}^{-ja}\right)t} \\
&= \lim_{t \rightarrow \infty} \mathrm{e}^{(\cos a -1)t}  \mathrm{e}^{jna}\mathrm{e}^{-jt \sin a} \\
&= 0.
\end{align*}
So, the cars rendezvous at the origin even though they are initially densely dispersed around the unit circle.}
\end{example}

\subsection{Convergence to the average starting point}

In the cyclic pursuit problem for a finite number of kinematic cars, the cars rendezvous at the average of their starting positions---see Example \ref{example_22}.
This turns out to be true also in the infinite chain serial pursuit problem provided there is an appropriate average initial position. The average of the points 
$\{ q_{m}(0), \dots, q_{m-N}(0)\}$ is
\[
\mathrm{avg} \{ q_{m-k}(0) \}_{k=0}^{N} = \frac{1}{N+1} \sum_{k=0}^{N} q_{m-k}(0).
\]
Our assumption will be that this average converges at the rate $1/\sqrt{N}$ as $N$ increases. That is, we will assume that 
there exists a number $\bar q$ such that
 \begin{equation} \label{eq_Cesaro}
 \mathrm{avg} \{ q_{m-k}(0) \}_{k=0}^{N}=\bar{q}+o(N^{-1/2}).
 \end{equation}
 This means that  for every constant $C$ there exists an integer $L$ such that if $N>L$ then
\[
\left|  \mathrm{avg} \{ q_{m-k}(0) \}_{k=0}^{N}-\bar{q} \right| \leq \frac{C}{\sqrt{N}}.
 \]

\begin{lemma}
 Assume that for some $m$ there exists a number $\bar q$ such that
(\ref{eq_Cesaro}) holds.
  Then $\lim_{t \rightarrow \infty} q_n(t) =\bar q$ for every $n$.
 \end{lemma}
 
 The proof uses a result of G.\ H.\ Hardy.
 Theorem 149 in \cite{Har91} is as follows: Let $\{ f_n\}_{n = 0}^{\infty}$ be a sequence of complex numbers such that there exists a constant $\bar{f}$ such that
\[
\mathrm{avg} \{f_0, \dots, f_n\} = \bar{f}+o(n^{-1/2}).
\]
Then 
\[
\lim_{t \rightarrow \infty} \mathrm{e}^{-t} \sum_{k=0}^{\infty} f_k \frac{t^k}{k!} = \bar{f}.
\]

\noindent
{\bf Proof of the Lemma} \ 
Assume (\ref{eq_Cesaro}) holds. Then it holds for all other values of $m$. For example
\begin{align*}
\mathrm{avg} \{ q_{m+1-k}(0) \}_{k=0}^{N} &= \mathrm{avg} \{ q_{m-k}(0) \}_{k=0}^{N} \\
& \ \ \ + \frac{1}{N+1} [q_{m+1}(0)-q_{m+1-N}(0)]\\
&= \bar{q}+o(N^{-1/2}) \\
& \ \ \ + \frac{1}{N+1} [q_{m+1}(0)-q_{m+1-N}(0)]\\
&= \bar{q}+o(N^{-1/2}).
\end{align*}
From Hardy's theorem, then, for every $m$
\[
\lim_{t \rightarrow \infty} \mathrm{e}^{-t} \sum_{k=0}^{\infty} q_{m-k}(0) \frac{t^k}{k!} = \bar{q}.
\]
From (\ref{eq_qn}), for every $m$
\[
\lim_{t \rightarrow \infty} q_{m}(t) = \bar{q}.
\]
This concludes the proof.

\subsection{An example of non-convergence}

It is more difficult to construct an example where $\lim_{t \rightarrow \infty} q_n(t)$ doesn't exist. We turn to such an example now.

In this example, the points $q_n(0)$ are either 0 or 1. Then existence of the limit $\lim_{t \rightarrow \infty} q_n(t)$ is related to some very interesting results of Diaconis and Stein on Tauberian theory 
 \cite{DiaSte78}, which we now briefly describe. Let $\mathbb{A}$ be an infinite subset of non-negative integers, for example the non-negative even integers, and consider the question of whether this limit exists:
 \begin{equation}
 \lim_{t \rightarrow \infty} \mathrm{e}^{-t} \sum_{k \in \mathbb{A}}\frac{t^{k}}{k!}.
 \label{lim1}
 \end{equation}
 Existence of the limit is a property of the set $\mathbb{A}$. Now let $S_n$ denote the number of heads that occur in $n$ tosses of a coin and consider the question of whether this limit exists:
 \begin{equation}
 \lim_{n \rightarrow \infty} \mathrm{Pr} (S_n \in \mathbb{A} ).
 \label{lim2}
 \end{equation}
 Finally, let $\mathit{card}$ denote cardinality and consider the question of whether this limit exists for every $\varepsilon >0$:
  \begin{equation}
 \lim_{n \rightarrow \infty} \frac{1}{\varepsilon \sqrt{n} } \ \mathrm{card} \{ k: k \in \mathbb{A}, n \leq k < n+\varepsilon \sqrt{n} \}.
 \label{lim3}
 \end{equation}
 Remarkably, the three limits are intimately related: If either exists, then so do the other two and they are all equal. This is Theorem 1 in   \cite{DiaSte78}.
 
 To get an example where (\ref{lim1}) fails, 
 by taking $n=m^2$ and $\varepsilon = 1$, it suffices to get an example where
 \[
  \frac{1}{ m } \ \mathrm{card} \{ k: k \in \mathbb{A}, m^2 \leq k < m^2 +m \}
  \]
  does not converge as $m \rightarrow \infty$. Defining
  \[
  \gamma(m)=\mathrm{card} \{ k: k \in \mathbb{A}, m^2 \leq k < m^2 +m \},
  \]
  it suffices to choose $\mathbb{A}$ such that $\gamma(m)=m-1$ for $m$ even and $\gamma(m)=0$ for $m$ odd. Returning to (\ref{eq_qn}), take $n=0$ and take the initial conditions $q_0(0), q_{-1}(0), \dots$ as follows. For $m\geq 0$ even, set $q_{-k}(0)=1$ for $m^2 \leq k < m^2 +m$, and for $m\geq 0$ odd, set $q_{-k}(0)=0$ for $m^2 \leq k < m^2 +m$. For other values of $k$, the value of $q_k(0)$ is irrelevant and could be set to 0.
Then $q_0(t)$ fails to converge as $t \rightarrow \infty$.


\subsection{Convergence of $q(t)$ on a subspace}

 Our final  result on this problem seems to be particularly interesting. 

\begin{theorem} \label{theorem4}
 The $\ell^\infty$-induced norm of $(U-I)\mathrm{e}^{(U-I)t}$ converges to 0 as $t \rightarrow \infty$.
Thus, for every $q(0) \in \ell^\infty$, $\dot{q}(t)$ converges to zero in $\ell^\infty$ as $t \rightarrow \infty$, and, moreover, if $q(0)$ belongs to $(U-I)\ell^\infty$, the image space of $U-I$ acting on $\ell^\infty$, then  ${q}(t)$ converges to zero in $\ell^\infty$ as $t \rightarrow \infty$.
\end{theorem}

\noindent
{\bf Proof} \
Given $q(0) \in \ell^{\infty}$, let $r(t)= (U-I)\mathrm{e}^{(U-I)t}q(0) $. To simplify layout, define
\[
\psi(k,t)=\frac{t^{k}}{k!} - \frac{t^{k+1}}{(k+1)!} .
\]
The $n^{\mathrm{th}}$ component of $r(t)$ is
\[
r_n(t)=\mathrm{e}^{-t}\left[-q_n(0) + \sum _{k=0}^{\infty}\psi(k,t) q_{n-(k+1)}(0)\right].
\]
Thus
\begin{align*}
\| r(t) \|_{\infty} &= 
\mathrm{e}^{-t}\sup_n \left|-q_n(0) + \sum _{k=0}^{\infty}\psi(k,t) q_{n-(k+1)}(0)\right|\\
& \leq  \mathrm{e}^{-t} \left[1 + \sum _{k=0}^{\infty}\left| \psi(k,t) \right| \right] \| q(0) \|_\infty .
\end{align*}
It therefore suffices to show that
\[
\lim_{t\rightarrow \infty} \mathrm{e}^{-t}\sum _{k=0}^{\infty}\left| \psi(k,t) \right|  = 0.
\]
It is elementary that for fixed $t$, the sequence $\{t^{k}/k!\}$ has a maximum at some $k_{0}(t)$, and that
this integer satisfies $k_{0}(t) \leq  t \leq k_{0}(t) + 1$. Also, for $k\leq k_{0}(t)$, the sequence is increasing
and for $k\geq k_{0}(t)$ the sequence is decreasing. Therefore,
\begin{align*}
\sum _{k=0}^{\infty}\left| \psi(k,t)\right| &=- \sum _{k=0}^{k_{0}(t) - 1}\psi(k,t)
 + \sum _{k=k_{0}(t)}^{\infty}\psi(k,t) \\\
&= 2\frac{t^{k_{0}(t)}}{k_{0}(t)! }- 1.
\end{align*}
Thus 
\[
\mathrm{e}^{-t}\sum _{k=0}^{\infty}\left| \psi(k,t)\right| \leq 2\mathrm{e}^{-t}\max_{k}\frac{t^{k}}{k!},
\]
which approaches zero as $t\rightarrow \infty$. This completes the proof.

\bigskip
Unfortunately, a simulation to illustrate the preceding result is not possible, because one cannot simulate an infinite number of kinematic points. The result implies that for every $q(0)$ in $\ell^\infty$
\[
\lim_{t \rightarrow \infty} \sup_n | q_{n-1}(t)-q_n(t) | =0.
\]
Intuitively, the  points $\{ q_n(t) \}$ cluster around a point, but that point may itself not be stationary. 

Theorem \ref{theorem4} raises the question of characterizing the image of the operator $U-I$.  
A bounded sequence $y$ belongs to $\mathrm{Im} (U-I)$ iff there exists a bounded $x$ such that $y=(I-U)x$, i.e.,
\[
x=(I+U+U^2+\cdots)y.
\]
Thus $\mathrm{Im} (U-I)$ equals the space of $y$ such that
$y$ and $  (I+U+U^2+\cdots)y$ are both bounded.

\section{A Symmetric Chain}

We turn now to kinematic cars where each is coupled to its two neighbours, not just one as in the preceding section.
Thus  the coupled equations
\begin{align*}
\dot{p}_n&=  (p_{n+1}-p_n)+(p_{n-1}-p_{n})
\end{align*}
or in vector form
\begin{equation}
\dot{p}=Ap, \ \ \ A=U+U^{-1}-2I.
\label{eq_example_diff}
\end{equation}
If the state space is $\ell^2$, then of course $U^{-1}=U^*$. The solution of (\ref{eq_example_diff}) is $p(t)=\mathrm{e}^{At}p(0)$.

\begin{theorem}\label{theorem1}
(With reference to (\ref{eq_example_diff}).)
\begin{enumerate}
\item
The $\ell^2$-induced norm of  $\mathrm{e}^{At}$ satisfies $\| \mathrm{e}^{At}\| = 1$ for all $t \geq 0$.
For every $p(0) \in \ell^2$, the $\ell^\infty$-norm of $p(t)$ converges to 0 as $t$ tends to $\infty$; 
in addition, $p(t)$ converges  to zero weakly, that is, the $\ell^2$ inner product $\langle p(t), y \rangle$ 
converges to zero as $t \rightarrow \infty$ for every $y \in \ell^2$.
\item
The $\ell^\infty$-induced norm of  $\mathrm{e}^{At}$ satisfies $\| \mathrm{e}^{At}\| = 1$ for all $t \geq 0$.
Also $\| A\mathrm{e}^{At}\| \rightarrow 0$ as $t\rightarrow \infty.$  Thus,  for every $q(0) \in \ell^\infty$, $\dot{q}(t)$ converges to zero in $\ell^\infty$ as $t \rightarrow \infty$, and, moreover, if $q(0)$ belongs to $A\ell^\infty$, the range space of $A$ acting on $\ell^\infty$, then  ${q}(t)$ converges to zero in $\ell^\infty$ as $t \rightarrow \infty$.
\end{enumerate}
\end{theorem}

\noindent
{\bf Proof} \
(1) The spectrum of $A$ is the real interval $[-4,0]$. The operator $FAF^{-1}$ is multiplication by 
\[
G(\mathrm{e}^{j \omega })=\mathrm{e}^{-j \omega }+\mathrm{e}^{j \omega }-2=2(cos \ \omega-1).
\]
The maximum real part of $G(\mathrm{e}^{j \omega })$ equals 0. Thus from (\ref{eq_norm}) $\| \mathrm{e}^{At}\| = 1$.

Take the spatial Fourier transform of the components of $p(t)$, holding $t$ fixed:
\[
P(\mathrm{e}^{j \omega},t)=\sum_n p_n(t) \mathrm{e}^{-j \omega n}.
\]
Then the function $\omega \mapsto P(\mathrm{e}^{j \omega},t)$, denoted for convenience by $P(t)$,  belongs to $\mathcal{L}^2(S^1)$. 
Taking Fourier transform of the differential equation gives
\[
\dot{P}(t)=(\mathrm{e}^{-j \omega }+\mathrm{e}^{j \omega }-2)P(t).
\]
Solve the equation:
\begin{align*}
P(t)& =\exp [ (\mathrm{e}^{-j \omega }+\mathrm{e}^{j \omega }-2)t] P(0) \\
& =\mathrm{e}^{ 2t ( \cos \omega -1)} P(0) .
\end{align*}
That is,
\begin{equation}
P(\mathrm{e}^{j\omega},t)=\mathrm{e}^{ 2t ( \cos \omega -1)} P(\mathrm{e}^{j\omega},0) .\label{eq_20}
\end{equation}
Recall that $\mathcal{L}^2(S^1)$ is a subspace of $\mathcal{L}^1(S^1)$. Using in turn the definition of the norm in $\mathcal{L}^1(S^1)$, equation (\ref{eq_20}), the Cauchy-Schwarz inequality, and  the modified Bessel function of the first kind $I_0$. we have
\begin{align*}
\| P(t) \|_1 & =\frac{1}{2 \pi} \int_{-\pi}^{\pi}   |P(\mathrm{e}^{j\omega},t) | d\omega \\
& =\frac{1}{2 \pi} \int_{-\pi}^{\pi}  \mathrm{e}^{ 2t ( \cos \omega -1)} |P(\mathrm{e}^{j\omega},0) | d\omega \\
& \leq   \left( \frac{1}{2 \pi} \int_{-\pi}^{\pi}  \mathrm{e}^{ 4t ( \cos \omega -1) }  d\omega \right)^{1/2} 
\| P(0) \|_2 \\
& = \left( \frac{1}{2 \pi}  \mathrm{e}^{-4t} I_0(4t) \right)^{1/2} 
\| P(0) \|_2.
\end{align*}
 But $\mathrm{e}^{-4t} I_0(4t)$ converges to 0 as $t \rightarrow \infty$ (the easiest way to see this is to graph the function), and therefore so does $\| P(t) \|_1$.  Since by definition
\[
p_n(t)=\frac{1}{2 \pi} \int_{-\pi}^{\pi}  P(\mathrm{e}^{j \omega },t) \mathrm{e}^{j \omega n} d\omega,
\]
so
\[
|p_n(t)| \leq \frac{1}{2 \pi} \int_{-\pi}^{\pi}  |P(\mathrm{e}^{j \omega },t) | d\omega =\| P(t) \|_1.
\]
Therefore $\| p(t) \|_\infty \leq \| P(t) \|_1$ and therefore $\| p(t) \|_\infty \rightarrow 0$ as $t \rightarrow \infty$.

Let $e_n$ denote the $n^{th}$ basis vector in $\ell^2$. We proved above that $p_n(t)$ converges to zero as $t \rightarrow \infty$. Thus
\[
\lim_{t \rightarrow \infty} \langle \mathrm{e}^{At} p(0),e_n \rangle =0.
\]
Therefore
\[
\lim_{t \rightarrow \infty} \langle \mathrm{e}^{At} p(0),y \rangle =0
\]
for every finite linear combination $y$ of $\{ e_n \}$. Also, $\| \mathrm{e}^{At} \| \leq 1$. Therefore by the uniform boundedness principle, 
\[
\lim_{t \rightarrow \infty} \langle \mathrm{e}^{At} p(0),y \rangle =0
\]
for every $y$ in $\ell^2$.

(2) Fix $t \geq 0$. The operator $U$ is an isometry. For every $x$ in $\ell^\infty$
\begin{align*}
\| \mathrm{e}^{Ut}x \|_\infty &= \left\| x+tUx +\frac{1}{2} t^2U^2 x + \cdots \right\|_\infty \\
& \leq \| x \|_\infty +t \| x \|_\infty +\frac{1}{2}t^2  \| x \|_\infty^2 + \cdots \\
& =\mathrm{e}^t \| x \|_\infty .
\end{align*}
And so $\| \mathrm{e}^{Ut} \| \leq \mathrm{e}^t.$ Likewise for $U^{-1}$. Thus
\[
\| \mathrm{e}^{At} \| \leq   \mathrm{e}^{-2t} \| \mathrm{e}^{Ut} \| \|  \mathrm{e}^{U^{-1}t} \| \leq  \mathrm{e}^{-2t} \mathrm{e}^{t} \mathrm{e}^{t}=1.
\]
To conclude equality, apply $\mathrm{e}^{At}$ to $\mathbf{1}$.

We next show that the kernel of $A=U+U^{-1}-2I$ as an operator on $\ell^\infty$ is the one-dimensional subspace spanned by $\mathbf{1}$.
Let $x$ be a bounded sequence in the kernel of $A$. Then
\[
x_{n}-x_{n-1}=x_{n+1}-x_n.
\]
In particular
\begin{align*}
x_{0}-x_{-1} &= x_1-x_0\\
x_1-x_0 & =  x_2-x_1 \\
& \mbox{ etc.}
\end{align*}
and so $x_{n+1}=x_n+(x_0-x_{-1})$ for all $n >0$. If $x_0-x_{-1} \neq 0$, then $x_n$ grows without bound as $n \rightarrow \infty$.
 Thus $x_0=x_{-1}$. Similarly $x_n=x_{n-1}$ for all $n$.
 
 Recall that $\| \mathrm{e}^{(U-I)t}\| \leq 1$ for all $t\geq 0$.
By symmetry the same holds for $\| \mathrm{e}^{(U^{-1} - I)t}\| $.
Thus 
\begin{align*}
 \| A\mathrm{e}^{At}\| &= \| [(U - I) + (U^{-1} - I)]\mathrm{e}^{(U-I)t}\mathrm{e}^{(U^{-1}-I)t}\| \\
&\leq \| \mathrm{e}^{(U^{-1} - I)t}(U-I)\mathrm{e}^{(U-I)t} \|  \\
& \ \ \ + \| \mathrm{e}^{(U-I)t}(U^{-1}-I)\mathrm{e}^{(U^{-1}-I)t}\|\\
& \leq \| (U-I)\mathrm{e}^{(U-I)t}\|  + \| (U^{-1}-I)\mathrm{e}^{U^{-1}-I)t}\| \\
& \rightarrow 0
\end{align*}
as $t\rightarrow \infty$. 
This completes the proof.

\bigskip
As our final contrast between the $\ell^2$ and $\ell^\infty$ cases, presented next is 
a chain that is unstable in $\ell^2$ but stable in $\ell^\infty$. The chain is not spatially invariant.

\begin{example}  \label{example_non_inv}
{\rm Define the operator  $B$ to repeat every component. Thus $B$ is defined by $y=Bx$, $y_{2n}=y_{2n+1}=x_n$.
As an operator on $\ell^2$, $\| B \| = \sqrt{2} $, $r_B=\sqrt{2}$, and $\sigma(B)=\{ \lambda : | \lambda | \leq \sqrt{2} \}$. As an operator on $\ell^\infty$, $\| B \| = 1 $, $r_B=1$, and $\sigma(B)=\{ \lambda : | \lambda | \leq 1 \}$. 

Let $a$ be a positive constant and consider the first-order system
\begin{align*}
\dot{p}_n &= u_n \\
u_{2n} &= -ap_{2n} +p_{n} \\
u_{2n+1} &= -ap_{2n+1} +p_{n}. \\
\end{align*}
Then
\[
\dot{p}=A p, \ \ \ A=-aI+B.
\]
This chain is not a simple mass-spring-dashpot system. The information flow structure is shown in Figure~\ref{fig_graph}. For example, for the component
\[
\dot{p}_5=-ap_5+p_2
\]
the graph shows an arrow from node 2 to node 5.
\begin{figure}[htbp] 
   \begin{center}
   \includegraphics[width=3in]{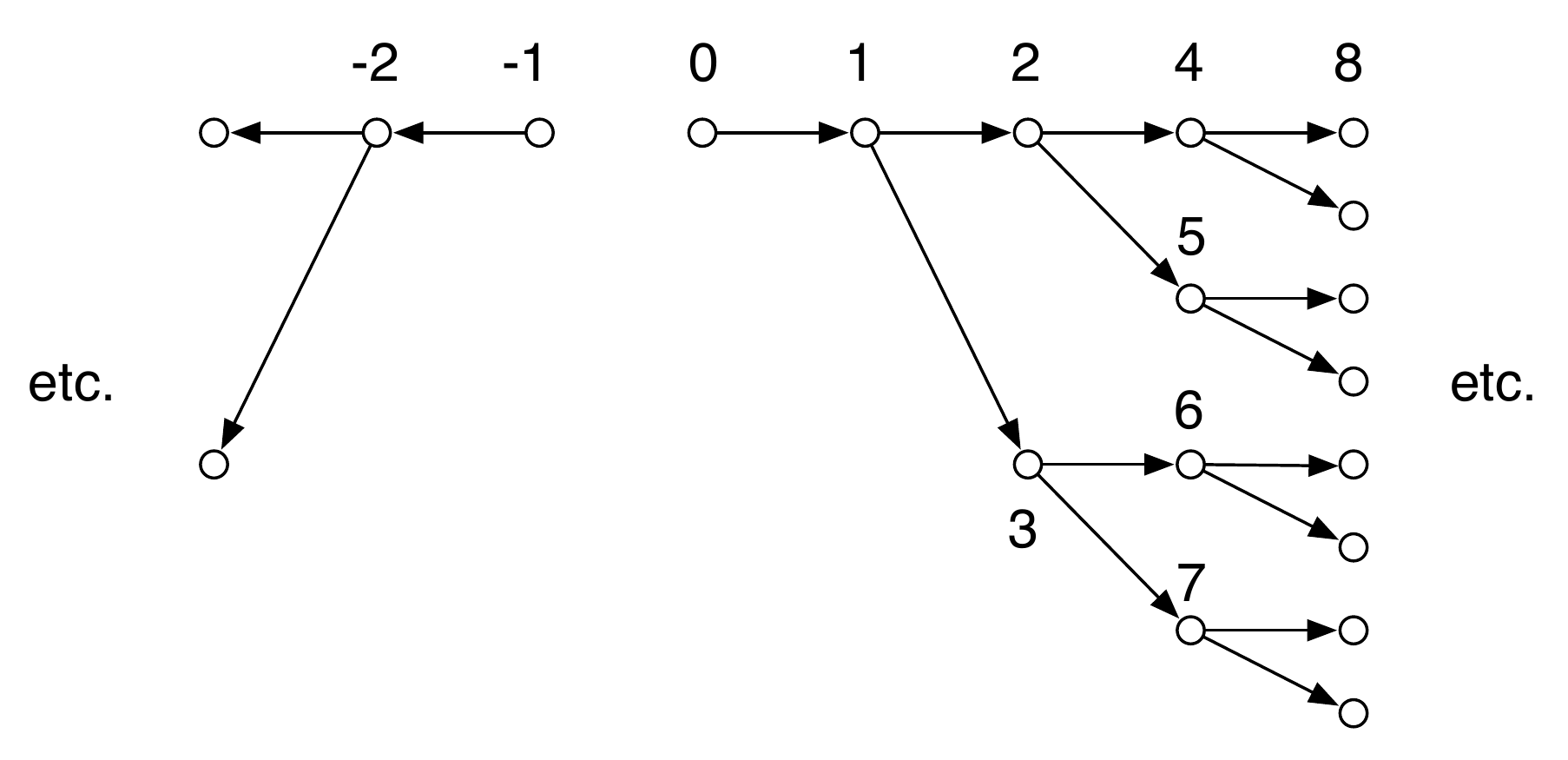} 
     \caption{Information flow in the example with spatial variation}
   \label{fig_graph}
     \end{center}
\end{figure}
The spectrum of $A$ equals
\[
\sigma(A) = -a + \sigma(B) .
\]
Select $a$ to lie in the interval $1 < a < \sqrt{2}$. Then, as an operator on $\ell^2$, $\sigma(A)$ has a nonempty intersection with the closed right half-plane, and so the origin $p=0$ {\bf is not} asymptotically stable; while on the other hand as an operator on $\ell^\infty$, $\sigma(A)$ is contained in the open left half-plane, and so the origin $p=0$ {\bf is} asymptotically stable. That is to say, there exists an initial state $p(0)$ in $\ell^2$ such that $\| p(t)\|_\infty$ converges to zero  but $\| p(t)\|_2$ diverges to $\infty$.}
\end{example}

\section{Literature review}

We now briefly review the literature. 
Chains are 1-dimensional lattices; lattices occur in physics problems. The theory of the propagation of waves, and in particular the application to determine the velocity of sound, is due to Newton and was published in 1687.
In Chapter III of \cite{Bri03}, Brillouin offers a mathematical treatment of wave propagation in a one-dimensional lattice of identical particles. However it is not mathematically rigorous. Reference \cite{KurHon95} is typical of the physics literature.
Kopell has a substantial oeuvre on chains of oscillators, for example \cite{Kop94}. In her work she does study the situation when the length of the chain goes to infinity.  However the boundary conditions are maintained.  We start with an infinite chain that therefore has no boundary conditions.

An early contribution to optimal control of an infinite chain of cars is that of Melzer and Kuo \cite{MelKuo71}. Their ``infinite object problem'' has the model
\[
\dot{x}(t)=Ax(t)+Bu(t),
\]
where for each $t$, $x(t)$ and $u(t)$ belong to $\ell^2$ and where $A$ and $B$ are spatially-invariant operators on $\ell^2$. The paper formulates a linear-quadratic optimal control problem with a cost function involving the time-domain $\mathcal{L}^2$-norm of 
\[
\langle x(t), Q u(t) \rangle + \langle u(t), Ru(t) \rangle ,
\]
 where $Q$ and $R$ are spatially-invariant operators on $\ell^2$. That is to say, the optimal control problem is formulated in the space $\mathcal{L}^2(\mathbb{R},\ell^2)$. The optimal control law takes the form $u=Fx$. The solution is derived via the Fourier transform. The work of Melzer and Kuo has been generalized and extended, most notably by 
 Bamieh et al.\ \cite{BamPagDah02},
D'Andrea and Dullerud \cite{DanDul03}, and
Motee and Jadbabaie \cite{MotJad08}. 
Curtain et al.\ \cite{CurIftZwa10} studied the LQR problem in the Hilbert space context, addressing the question of truncating the infinite chain. 

Our paper is related to that of Swaroop and Hedrick on string stability \cite{SwaHed96}. The system in that paper is a semi-infinite chain, that is, the cars are numbered $0,1,2,\dots$ and car 0 is therefore a boundary, its dynamics being independent of all others. This model is appropriate for a platoon (or convoy) with a leader. By contrast, in the other references and in our paper there is no boundary car. On the other hand, reference  \cite{SwaHed96} is the only reference we found that proposes $\ell^\infty$ for the state space.

Other recent papers are \cite{BamVou05},  \cite{Cur09},  \cite{HuiBer09}, and \cite{jovbam05}.

\section{Conclusion}

An infinite chain of vehicles obviously doesn't occur in reality. It is intended to be relevant to the case of a finite but very large chain. We have argued that the $\ell^2$-framework is not the right one, because an infinite chain does not behave like a finite but large one. For example, formulating the rendezvous problem in $\ell^2$ results in a rendezvous at the origin, whereas a finite chain would not do so. The $\ell^\infty$ formulation seems more appropriate in our opinion. However, the problems are harder. Designing controllers that are $\ell^\infty$-optimal is an open problem.

\subsection*{Acknowledgements}

Thanks to Harry Dym for drawing our attention to \cite{DalKre70},  to Ronen Peretz for introducing us to and discussions on \cite{DiaSte78}, to Joyce Poon for introducing us to the physics literature, and to Geir Dullerud and Bassam Bamieh for comments about the existing literature.

\bibliographystyle{plain}
\bibliography{bf}

\end{document}